\documentstyle [11pt,graphicx,amstex,epsfig,amscd,amssymb]{article}
\textwidth
6.2in

\def\om{\omega}

\def\phi{\varphi}
\def\ro{\rho}
\def\g{\mathfrak g}

\def\ti{\times}

\def\be{\begin{equation}}
\def\ee{\end{equation}}
\def\bear{\begin{eqnarray}}
\def\eear{\end{eqnarray}}
\def\best{\begin{eqnarray*}}
\def\eest{\end{eqnarray*}}
\def\pf{{ {\emph{Proof}}}: }
\def\ni{\noindent}

\newtheorem{theorem}{Theorem}[section]
\newtheorem{prop}[theorem]{Proposition}
\newtheorem{proposition}[theorem]{Proposition}

\newtheorem{cor}[theorem]{Corollary}

\newtheorem{defn}[theorem]{Definition}

\newtheorem{remark}[theorem]{Remark}

\newtheorem{example}[theorem]{Example}

\renewcommand{\Re}{\operatorname{Re}}

\renewcommand{\Im}{\operatorname{Im}}

\def\ra{\rightarrow}

\def\r#1{\right#1}
\def\l#1{\left#1}

\def\R{{ \Bbb R}}
\def\C{{ \Bbb C}}

\def\ts3{{T^*S^3}}
\def\tsn{{T^*S^n}}

\title{\bf Cohomogeneity One Special Lagrangian Submanifolds in the Deformed Conifold}

\author{Marianty Ionel and Maung Min-Oo}

\date{}

\begin{document}


\maketitle

\vskip.1in

\noindent {\bf Abstract. } In this paper we describe the
cohomogeneity one special Lagrangian 3-folds in the cotangent
bundle of the $3$-sphere, also known in the physics literature as
a deformed conifold.  Our main result gives a global foliation of
the deformed conifold by $T^2$-invariant special Lagrangian
$3$-folds, where the generic leaf is topologically $T^2 \times
\R$. In the limit, these special Lagrangians asymptotically
approach a special Lagrangian cone on a torus in the conifold.
Using moment map techniques we also recover the family of
$SO(n)$-invariant special Lagrangian $n$-folds in $\tsn$ obtained
by H. Anciaux in \cite{an}.
\\

\\


\section{Introduction}
\bigskip

Beginning with the seminal paper by R.Harvey and H.B. Lawson
\cite{hl} on calibrated geometry, there has been extensive
research in the mathematics literature on special Lagrangian and
other calibrated submanifolds. Recently, a lot of progress has
been done in constructing special Lagrangian submanifolds using
various techniques. To give some examples, D. Joyce used the
method of ruled submanifolds, integrable systems and evolution of
quadrics in \cite{dj3,dj4,dj5} to construct
 explicit examples of special Lagrangian $m$-folds in $\C^m$, M. Haskins exhibited
examples of special Lagrangian cones in $\C^3$ \cite{ha}, etc.
Although the main emphasis has been on examples in flat space,
such as $\C^n$, there has also been progress in studying the
special Lagrangian submanifolds in non-flat Calabi-Yau manifolds.
R. Schoen and J. Wolfson used variational methods to construct
special Lagrangians in Calabi-Yau manifolds in \cite{sw}. H.
Anciaux found new $SO(n)$-invariant examples in $\tsn$ \cite{an},
equipped with the Ricci-flat Stenzel metric. See also \cite{ikm}
and \cite{km}, where a different method is used to construct
calibrated submanifolds.

In this paper, we will find new examples of cohomogeneity one
special Lagrangians in the Calabi-Yau manifold $\ts3$, which is
known in the physics literature as the deformed conifold. The main
result is that we exhibit a foliation of $\ts3$ by special
Lagrangians, where the generic leaf is $T^2 \times \R$ and the
$T^2$ is an orbit under the maximal torus of the isometry group
$SO(4)$. We will show that asymptotically these special
Lagrangians converge to a special Lagrangian cone in the conifold,
which is the cone on $S^2 \times S^3$, with the singular
Calabi-Yau metric \cite{co}. We will also recover some of
Anciaux's results about $SO(3)$-invariant special Lagrangians
using different techniques. Our methods use the moment map of the
$SO(4)$ action and are similar to Joyce's methods \cite{dj2} for
the flat case.

The Calabi-Yau metric on $\ts3$, the deformed conifold, was first
explicitly described by Candelas and de la Ossa \cite{co} in 1990,
before Stenzel (independently) discovered it for any $T^*S^n$, and
more generally, for any cotangent bundles of rank one symmetric
spaces. The case $T^*S^2$ is already very interesting and the
Ricci-flat hyperk{\"a}hler metric on this manifold was first
discovered by Eguchi and Hanson \cite{eh} in 1978. This
4-dimensional metric is the basic model of a number of explicit
special holonomy metrics discovered by physicists.

Special Lagrangian submanifolds in Calabi-Yau spaces play an
important role in string theory, particularly in mirror symmetry.
Since Calabi-Yau metrics are difficult to describe explicitly on
compact manifolds, it is a good strategy to understand the local
non-compact models first. Moreover, it is important to study how
singularities develop, since they provide a transition between
different topological types in the moduli space of Calabi-Yau
structures. In dimension $6$, the basic example is that of a
conifold transition, where a Calabi-Yau metric on a smooth
manifold diffeomorphic to $T^*S^3$, called the deformed conifold,
degenerates into a singular cone metric on $S^2 \times S^3$ and
then is resolved (by a blow-up) into another Calabi-Yau metric on
the total space of a $\C^2$-bundle $({\cal O}(-1) + {\cal O}(-1))$
on $S^2$, which is known as the {\it small resolution} of the
conifold and is the mirror of the deformed conifold. This
transition and the metrics involved are described explicitly by
Candelas and de la Ossa \cite{co}.

We now give a brief description of the contents of this paper. In section 2, we review
the basic facts about Calabi-Yau manifolds and special Lagrangians, including a list of well-known
 relevant examples. In sections 3 we describe the Calabi-Yau structure on the deformed conifold and
 the Stenzel metric on $\tsn$. We define the associated moment map and discuss some of its
  basic properties in section 4. In section 5 we prove that the only
 homogenous example of special Lagrangians $3$-fold in $\ts3$ is the zero section and in
 section 6 we compute our main examples
 of cohomogeneity one special Lagrangians including a description of the
 asymptotics. We conclude with some remarks and open questions in
 the last section.

\section{Special Lagrangian Geometry}

\bigskip

Special Lagrangian submanifolds are a special class of minimal
submanifolds in Calabi-Yau spaces and were introduced by Harvey
and Lawson in their seminal paper \cite{hl} using the notion of a
calibration. We begin by quickly reviewing the basic definitions
and setting up notations. For details, see \cite{dj1}.
\subsection{Basic Definitions}

\begin{defn}
A Calabi-Yau $n$-fold $(M,J,\omega,\Omega)$ is a K{\"a}hler
$n$-dimensional manifold $(M,J,\omega)$ with Ricci-flat K{\"a}hler
metric $g$ and a nonzero holomorphic section $\Omega$ which
trivializes the canonical bundle $K_M$.
\end{defn}

Since the metric $g$ is Ricci-flat, $\Omega$ is a parallel tensor
with respect to the Levi-Civita connection $\nabla^g$ \cite{dj1}.
By rescaling $\Omega$, we can take it to be the holomorphic
$(n,0)$-form that satisfies:
\begin{equation}
{\frac{\om^n}{n!}}={(-1)^{\frac{n(n-1)}{2}}{\biggl(\frac{i}{2}\biggl)}^n\Omega\wedge\bar
\Omega}, \label{holo}
\end{equation}
where $\om$ is the K{\"a}hler form of $g$. The form $\Omega$ is
called the {\it holomorphic volume form} of the Calabi-Yau
manifold $M$.

Let $\phi$ be a closed $p$-form on manifold $M$. We say $\phi$ is
a calibrating form on $M$ if $$\phi_{\mid V}\leq vol_V$$ for any
oriented $p$-plane $V\subset T_xM, \forall x\in M$. A submanifold
$N$ of $M$ is called {\it calibrated by $\phi$} if $\phi\mid_{T_x
N}= vol_{T_x N}, \forall x\in N$.
\smallskip

\noindent{\it Remark:} The constant factor in $(\ref{holo})$ is
chosen so that $\Re \Omega$ becomes a calibration on $M$.

\begin{defn}
Let $(M,J,\omega,\Omega)$ be a Calabi-Yau $n$-fold and $L\subset
M$ a real oriented $n$-dimensional submanifold of $M$. Then $L$ is
called a {\it special Lagrangian submanifold} of $M$ if it is
calibrated by $\Re \Omega$.
\end{defn}
\ni More generally, $L$ is said to be a special Lagrangian
submanifold with phase $e^{i\theta}$ if it is calibrated by $
\Re(e^{-i \theta} \Omega)$, where $\theta$ is a constant.

\medskip
An alternative and very useful description of the special
Lagrangian submanifolds found in \cite{hl}  is the following:

\begin{prop}
Let $(M,J,\omega,\Omega)$ be an $n$-dimensional Calabi-Yau
manifold and $L\subset M$ a real $n$-dimensional submanifold of
$M$. Then $L$ is called a {\it special Lagrangian submanifold} of
$M$ if $ \om\mid_L\equiv 0$ and $\Im\Omega\mid_L\equiv 0$.
\end{prop}

\noindent {\it Remark:} The condition $ \om\mid_L\equiv 0$ says
that $L$ is a Lagrangian submanifold. Therefore, special
Lagrangian submanifolds are Lagrangian with the extra condition
$\Im\Omega\mid_L\equiv 0$.
\smallskip

 \noindent The simplest example of a Calabi-Yau manifold
is $\C^n$, with coordinates $(z_1,\dots,z_n)$, endowed with the
flat metric $<,>$, the K{\"a}hler form $\omega_0$ and the
holomorphic volume form $\Omega_0$, where:

\begin{align}
\label{str1} <z,w> & =\Re \, \sum_{i=1}^n z_i\bar{w_i}  \\
\label{str2} \om_0 & =\frac{i}{2}(dz_1\wedge d\bar z_1+\dots
+dz_n\wedge d\bar z_n)\\
\label{str3} \Omega_0 & = dz_1\wedge\dots\wedge dz_n
\end{align}

\noindent $\R^n$ is a trivial example of a special Lagrangian
submanifold in $\C^n$.

One important property of the special Lagrangian submanifolds is
that they are absolutely area minimizing in their homology class,
so they are in particular {\it minimal} submanifolds (see
\cite{hl}).

The graph of $F:\R^n \to \R^n$ is a Lagrangian $n$-fold in
$\R^{2n}$ if and only if $F=\nabla f$, for some function $f:\R^n
\to \R$. The graph is special Lagrangian if and only if the
function $f$ satisfies the differential equation:
$$\Im det_{\C}(I+i \mbox {Hess }f)=0 \mbox{  on } \C^n$$

\noindent This is a nonlinear elliptic P.D.E. and it is difficult
to solve in general. One idea, initiated by Harvey and Lawson, is
to look for solutions invariant under certain group actions of
$\C^n$.  In their paper \cite{hl}, Harvey and Lawson produced many
examples of symmetric special Lagrangian. We will review some of
their examples relevant to this paper in the next section.

\subsection{Cohomogeneity one examples in $\C^n$}\label{sectex}

Let $L$ be a special Lagrangian submanifold of $\C^n$. The group
of automorphisms preserving the Calabi-Yau structure on $\C^n$ is
$G=SU(n)\ltimes \C^n$. The symmetry group of $L$ is defined to be
the Lie subgroup of $SU(n)\ltimes \C^n$ that acts on $\C^n$
leaving $L$ fixed. In general, the easiest special Lagrangian
submanifolds to construct are those with large symmetry groups. It
can be shown that all homogeneous special Lagrangian submanifolds
of $\C^n$ are conjugate under $SU(n)\ltimes \C^n$ to the standard
real $n$-plane $\R^n\subset \C^n$ \cite{dj1}. The next most
symmetric special Lagrangians are those of cohomogeneity one, i.e.
the ones whose orbits of the symmetry group are of codimension one
in $L$. We will now recall two important examples of cohomogeneity
one special Lagrangian submanifolds found by Harvey and Lawson in
\cite{hl}.

\vskip.2cm

\noindent{\bf Example 1:} The following are special Lagrangian
submanifolds of $\C^n$, invariant under the action of $SO(n)
\subset SU(n)$:
\begin{equation}
\label{eqHL}
 L_c= \{ \lambda \, u  \,\,|\, \, u\in S^{n-1}\subset\R^n, \,\, \lambda \in \C,
\ \Im(\lambda^n)=c \, \},
\end{equation}
\ni where $c\in \R$ is a constant. The variety $L_0$ is an union
of $n$ special Lagrangian $n$-planes and when $c\not = 0$, each
component of $L_c$ is diffeomorphic to $\R\ti S^{n-1}$ and it is
asymptotic to $L_0$.

\medskip
\noindent {\bf Example 2:} The special Lagrangian submanifolds in
$\C^n$, invariant under the action of the maximal torus
$T^{n-1}=\{\mbox{ diag}(e^{i\theta_1}, \dots, e^{i\theta_n})|
\theta_1+ \dots \theta_n=0\}$ of $SU(n)$ are given by the
equations:

\begin{align*}
&|z_1|^2-|z_j|^2=c_j, \,\,j=2,3,\dots,n \mbox{ \quad    and}\\
& \Re (z_1 z_2 \dots z_n)=c_1,\, \mbox {  if $n$ even}\\
& \Im (z_1 z_2 \dots z_n)=c_1,\, \mbox {  if $n$ odd}
\end{align*}
where $c_1,c_2,\dots c_n$ are any real constants.
\smallskip

\ni {\it Remark:} These special Lagrangians are topologically
$T^{n-1} \times \R\,$ and they are asymptotic to the union of two
cones on flat tori $T^{n-1}$ in $S^{2n-1}$ obtained by making all
constants equal to $0$.

\smallskip

Using moment map techniques, D. Joyce \cite{dj1} made a systematic
study of cohomogeneity one special Lagrangian folds in $\C^3$ and
he found the following new example:

\smallskip
\noindent {\bf Example 3:} The submanifold given by
$$L_{a,b,c}=\{(z_1,z_2,z_3)\in\C^3: |z_1|^2-|z_2|^2=a, \
\Re(z_1z_2)=b, \ \Im(z_3)=c\}$$ where $a,b,c$ are any real
constants is a special Lagrangian in $\C^3$ and it is invariant
under the subgroup $U(1)\ltimes \R$ of $SU(3)\ltimes \C^3$, acting
by
$$(e^{i\theta},x).(z_1,z_2,z_3)\to (e^{i\theta} z_1,e^{-i\theta}
z_2,x+z_3)$$ for $\theta\in[0,2\pi), x\in\R$.

D. Joyce in \cite{dj2} proved that the above examples are all the
cohomogeneity one special Lagrangian submanifolds in $\C^3$:

\begin{theorem}

Every cohomogeneity one special Lagrangian $3$-fold in $\C^3$ is
conjugate under $SU(3)\ltimes \C^3$ to a subset of $\R^3$, one of
the $3$-folds of Example 1 and Example 2 with $n=3$ or one of
those in Example 3.

\end{theorem}

\section{$\tsn$ as a Calabi-Yau manifold}

In 1995, M. Stenzel \cite{st} discovered that the cotangent bundle
of the sphere, $T^*S^n$, and more generally the cotangent bundle
of symmetric spaces of rank one, can be endowed with a Ricci-flat
metric making them into Calabi-Yau manifolds. As mentioned in the
introduction, the lower dimensional cases were discovered in the
physics literature. The case $n=2$ is the Eguchi-Hanson metric
\cite{eh} and the case $n=3$ (the deformed conifold) is due to
Candelas and de la Ossa \cite{co}.

\noindent In what follows we will first describe the cotangent
bundle of the sphere as a complex affine quadric following
Sz{\"o}ke \cite{rs} and define the K{\"a}hler potential of the
Stenzel metric. For more details on this material, the reader
should consult \cite{st}.

\subsection{The Stenzel metric}\label{cy}
 Let $\tsn=\{(x,\xi)\in\R^{n+1}\times \R^{n+1}, \ |x|=1, \ x\cdot
\xi=0\}$ be the cotangent bundle of the $n$-sphere (which we have
identified with the tangent bundle). The group $SO(n+1,\R)$ acts
with cohomogeneity one on $\tsn$, the generic orbit being
$SO(n+1)/SO(n-1)$ (where $|\xi|$ is constant). According to Sz{\"
o}ke \cite{rs}, one can identify $\tsn$ with the affine quadric
$$Q^n=\{z=(z_0,z_1,\dots,z_{n})\in\C^{n+1}\big | \ \sum_{i=0}^{n}
z_i^2=1\}$$ using the diffeomorphism $h:\tsn\to Q^n$ given by:
\begin{equation}\label{h}
(x,\xi)\ra z=x\cosh(|\xi|)+ i \ \frac{\sinh(|\xi|)}{|\xi|}\xi
\end{equation}
This diffeomorphism is equivariant with respect to the action of
$SO(n+1,\R)$ on $\tsn$ and the natural action of $SO(n+1,\C)$ on
$Q^n$. The complex structure on the cotangent bundle of the
$n$-sphere is obtained by pulling back the complex structure of
the affine quadric under the map $h$. On the complex quadric,
there exists a Ricci-flat metric whose corresponding symplectic
form is the Stenzel form given by:
$$\omega_{St}=i\partial\bar{\partial}
u(r^2)=i\sum_{j=0}^{n}\sum_{k=0}^{n} \frac{\partial^2}{\partial
z_j
\partial\bar{z_k}} u(r^2) dz_j\wedge d \bar{z_k}$$
where $r^2=|z|^2=\sum_{j=0}^n z_j\bar{z_j}=\cosh(2|\xi|)$ and
$u(r^2)$ is a smooth real function satisfying the following
differential equation:
\begin{equation}\label{ode}
\frac{d}{d\tau}(u'(\tau))^n= c\,n\,(\sinh \tau)^{n-1}
\end{equation}

\ni where $\tau =\cosh^{-1}(r^2)$ and $c$ is a positive constant
(see \cite{st}). In dimension $n=2$, there is an explicit formula
for the potential function: $u(r^2)=\sqrt{1+r^2}$. In dimension
$n=3$, the derivative of the potential function is given by the
relation:
\begin{equation}\label{diffu}
u'(\tau)^3=\frac{3c}{2}\l(\frac{\sinh(2\tau)}{2}-\tau\r)
\end{equation}
which using the initial condition $u'(0)=0$ integrates to the more
complicated formula:
\begin{equation}\label{diff3}
u(r^2)=\int_0^{\cosh^{-1}(r^2)}\l[\frac{3c}{2}\l(\frac{\sinh(2\sigma)-\sigma}{2}\r)\r]^{\frac{1}{3}}d\sigma
\end{equation}

\smallskip
\noindent The form $\omega_{St}=i\partial\bar{\partial}u$ is exact
on $T^*S^n$ and $\omega_{St}=d\alpha_{St}$, where
$\alpha_{St}=-\Im (\bar{\partial} u)$. The form $\alpha_{St}$ is
related to the Liouville form $\alpha_0(v)=\frac{1}{2}<v,Jz>$ on
$\C^{n+1}$ by $\alpha_{St}=u'(|z|^2)\alpha_0$. Therefore, it
follows that the $1$-form $\alpha_{St}$ has the expression:
\begin{equation}
\alpha_{St}(v)=\frac{1}{2}u'(|z|^2)\omega_0(z,v)=\frac{1}{2}u'(|z|^2)<v,Jz>,
\ v\in T_z Q, z \in Q\label{alfa}
\end{equation}

\ni where $<,>$ and $\omega_0$ are respectively the flat Euclidean
metric and the K{\"a}hler form of $\C^{n+1}$. It is well known
that on the complex space $\C^{n+1}$, $\omega_0(v,w)=<Jv,w>$,
where $J$ is the complex structure on $\C^{n+1}$.

Differentiating expression (\ref{alfa}), we calculate:
\begin{align*}
\omega_{St}(v,w)&=d\alpha_{St}(v,w)=v(\alpha_{St}(w))-w(\alpha_{St}(v))-\alpha_{St}([v,w])\\
&=v\l(\frac{1}{2}u'(|z|^2)<w,Jz>\r)-w\l(\frac{1}{2}u'(|z|^2)<v,Jz>\r)-\alpha_{St}([v,w])\\
&=\frac{1}{2}\{u''(|z|^2)v(|z|^2)<w,Jz>+u'(|z|^2)v(<w,Jz>)- u''(|z|^2)w(|z|^2)<v,Jz>\\
&-u'(|z|^2)w(<v,Jz>)-u'(|z|^2)<[v,w],Jz>\}\\
&=u'(|z|^2)<w,Jv>+u''(|z|^2)(<v,z><w,Jz>-<w,z><v,Jz>)\\
&=u'(|z|^2)\omega_0(v,w)+u''(|z|^2)(<v,z>\omega_0(z,w)-<w,z>\omega_0(z,v))
\end{align*}

\ni In the above calculation we used that $v(|z|^2)=2<v,z>$,
$\nabla_v w-\nabla_w v=[v,w]$ and the fact that $<v,Jw>=-<w,Jv>$.

 \ni Therefore, the K{\"a}hler
form of the Stenzel metric at a point $z$ on the quadric $Q$ is
given by:
\begin{align}
\omega_{St} (v,w) =  u'(|z|^2)\,\omega_0(v,w) +
u''(|z|^2)\Big(<w,z>\omega_0(v,z) - <v,z>\omega_0(w,z)  \Big) , \
v,w\in T_z Q \label{St2}
\end{align}

\ni Formulas (\ref{alfa}) and (\ref{St2}) will prove to be
important when we compute the moment maps for group actions on the
quadric.

\ni On the quadric $Q$ define the holomorphic $(n,0)$-form
$\Omega_{St}$ by the relation:
$$\frac{1}{2}d(z_0^2+z_1^2+\dots+z_n^2-1)\wedge \Omega_{St}=\Omega_0$$
where $\Omega_0 = d z_0 \wedge d z_1 \wedge \dots \wedge d z_n$ is
the holomorphic volume form of $\C^{n+1}$. Therefore:
\begin{align}
\Omega_{St}(v_1,v_2,\dots, v_n)=&\Omega_0 \l( z,
v_1,\dots,v_n\r),\ v_1,\dots, v_n\in T_zQ, \ z\in Q \label{holom}
\end{align}

\ni The quadric $Q^n$ with this structure becomes a Calabi-Yau
manifold since equation (\ref{holo}) holds for $\omega_{St}$ and
the corresponding holomorphic $n$-form $\Omega_{St}$, up to a
multiplicative constant (see also \cite{an}).

\subsection{The conifold in dimension 3}\label{sect}
Let $Q_0$ be the quadric in $\C^{4}$ defined by the equation:
$$\sum_{i=0}^{3}z_i^2=0$$ This quadric, called the conifold, is singular at the origin and represents
a cone on $T_1(S^3)\cong S^2\times S^3$. The complex structure of
the conifold is given by the embedding
$$h_0:T_1(S^3)=\l\{(x,\xi)\in\R^4\times\R^4  | \ |x|=|\xi|=1, \ x\cdot
\xi=0\r\}\to Q_0, \ z=h_0(x,\xi)=x+i\xi$$ Deforming the conifold
equation to $\sum_{i=0}^{3}z_i^2=\epsilon^2$, with $\epsilon$ a
positive constant, yields the complex quadric $Q_{\epsilon}$ in
dimension 3, where $\epsilon$ is the radius of the zero section
$S^3$. This is equivalent to replacing the tip of the conifold by
an $S^3$( see \cite{co}) and as $\epsilon\to 0$, this sphere
collapses into the singular point of $Q_0$. In the physical
literature, the complex quadric $Q_\epsilon$ is also known as a
{\it deformed conifold}.

\ni Candelas and de la Ossa \cite {co} showed that the conifold
$Q_0$ admits a Ricci-flat metric $g_{cone}$ with K{\"a}hler
potential $u_{cone}(r^2)=\frac{3}{2}r^{\frac{4}{3}}$. The
holomorphic $(n,0)$-form $\Omega_{cone}$ defined by:
\begin{equation}
\frac{1}{2}d(z_0^2+z_1^2+\dots+z_3^2)\wedge \Omega_{cone}=d z_0
\wedge d z_1 \wedge \dots \wedge d z_3 \label{holcone}
\end{equation}
makes the conifold into a singular Calabi-Yau manifold.

We note that the above relation can be used to compute
$\Omega_{cone}$ as:
\begin{equation}\label{hol}
\Omega_{cone}(v_1,v_2,v_3)=\frac{1}{|z|^2}(dz_0\wedge dz_1\wedge
dz_2\wedge dz_3)(\bar{z},v_1,v_2,v_3)
\end{equation}
where $v_1,v_2,v_3\in T_z Q_0$,  $z\in Q_0$ and
$\bar{z}=\bar{z_0}\frac{\partial}{\partial z_0}+
\bar{z_1}\frac{\partial}{\partial
z_1}+\bar{z_2}\frac{\partial}{\partial
z_2}+\bar{z_3}\frac{\partial} {\partial z_3}$.

\ni We make the remark that we can not use the position vector $z$
anymore on the cone in this case to calculate $\Omega_{cone}$ in
terms of $\Omega_0$, since
$d(z_0^2+z_1^2+z_2^2+z_3^2)\wedge\Omega_{cone}(z,v_1,v_2,v_3)=0$.
Instead, we use the vector $\bar{z}$ which is normal to the cone
and the fact that
$\frac{1}{2}d(z_0^2+z_1^2+\dots+z_3^2)(\bar{z})=|z|^2$.

\section{Moment Maps and Special Lagrangians with Symmetry}

The group $SO(n+1,\R)$ acts on the cotangent bundle of the sphere
\begin{equation}\label{cot}
 T^*S^n=\{(x,\xi)\in\R^{n+1}\times \R^{n+1}, \ |x|=1, \
x\cdot \xi=0\}
\end{equation}
 with cohomogeneity one. The action is transitive on
the sets $|\xi|= \rho =\mbox{constant}$ and is given by:
\begin{equation}\label {act}
g.(x,\xi)=(gx,g\xi), \ g\in SO(n+1), \ (x,\xi)\in T^*S^n
\end{equation}

\ni We would like to find examples of special Lagrangian
submanifolds in the deformed conifold $T^*S^3$, with large
symmetry group. Our method will use moment map techniques. In what
follows we will compute the moment map of a group action on the
cotangent bundle of the sphere and review some results of Joyce
\cite{dj2} about finding $G$-invariant special Lagrangian
submanifolds in $\C^n$ by using the moment map. These results are
quite general and easily extend from the case of the flat $\C^n$
to general symplectic manifolds with Hamiltonian actions and in
particular to the Calabi-Yau manifold $T^*S^n$.

Let $Q \cong \tsn$ be endowed with the Calabi-Yau structure
described in section \ref{cy}. As we have seen, the group of
automorphisms of $\tsn$ preserving the Calabi-Yau structure is
$SO(n+1,\R)\subset SO(n+1,\C)$. Let $G$ be a Lie subgroup of
$SO(n+1,\R)$, with Lie algebra $\mathfrak g $. Let $A
\in{\mathfrak g} \subset {\mathfrak o}(n+1)$. Then the induced
vector field on $Q$ is given by: $z \mapsto X_A(z) = Az $ with
flow $z \mapsto e^{tA}z$ for $z \in Q$. When the symplectic form
is exact, like in the case of the Stenzel form
$\omega_{St}=d\alpha_{St}$, the moment map $\mu: T^*S^n\to
{\mathfrak g}^* $ of the $G$-action can be computed easily as
follows:

\ni Both forms $\omega_{St}$ and $\alpha_{St}$ are invariant under
the flow of $X_A$, that is $${\mathcal L}_{X_A}\omega_{St} =
{\mathcal L}_{X_A}\alpha_{St} = 0$$

\ni for any $A \in {\mathfrak g} \subset {\mathfrak o}(n+1)$.

\ni From Cartan's formula,  ${\mathcal L}_{X_A}\alpha_{St}=d(X_A
\, \lrcorner \, \alpha_{St})+ X_A \, \lrcorner \, d\alpha_{St}$
and using $d\alpha_{St}= \omega_{St}$ , we see that
$$ - X_A \,\lrcorner \, \omega_{St} = d(\alpha_{St}(X_A)) $$

\ni This shows that the action of $G$ on $\tsn$ is hamiltonian,
with moment map $\mu(z)$ given by $A \mapsto \alpha_{St}(A)$ which
can be rewritten as:
$$
<\mu(z), A>= \mu_A(z) = \alpha_{St}(X_A(z))$$ where $z\in Q
\subset \C^{n+1}$, $A \in  {\mathfrak g} \subset {\mathfrak
o}(n+1)$ and $<,>$ is the pairing between $\mathfrak g$ and
$\mathfrak g^*$. Using formula (\ref{alfa}) for $\alpha_{St}$, we
have thus proved the following proposition.

\smallskip
\begin{proposition}
Let $G\subset SO(n+1)$ be a connected Lie group. The moment map of
the $G$-action on the complex quadric $Q$ is given by:
\begin{equation}\label{action}
\mu_A:Q^n\to \R, \quad
\mu_A(z)=\alpha_{St}(Az)=\frac{1}{2}u'(|z|^2)<Az,Jz>, \ \forall
A\in {\mathfrak g}
\end{equation}
where $<,>$ is the Euclidean metric on $\C^{n+1}$ and the function
$u$ satisfies equation $(\ref{ode})$.
\end{proposition}

\medskip
\ni We define the center $Z(\mathfrak g)$ of $\mathfrak g$ to be
the subspace of $\g$ fixed by the coadjoint action of $G$. Since
the moment map is equivariant, a level set of the moment map
$\mu^{-1}(c)$ for $c\in\mathfrak g^*$ is $G$-invariant if and only
if $c\in Z(\g ^*)$. \ni In order to find examples of special
Lagrangian submanifolds in $Q$, we will use the following
proposition which is true for moment maps on any symplectic
manifold, in particular on $Q$.

\begin{proposition}\label{lagr}
\label{mm} Let $G\subset SO(n+1)$ be a connected Lie subgroup with
Lie algebra $\g$ and moment map $\mu:\tsn\to \g^*$ and $\mathcal
O$ an orbit of $G$ in $\tsn$. Then the orbit is isotropic, i.e.
$\omega_{St}|_{\mathcal O}\equiv 0$ if and only if ${\mathcal
O}\subseteq \mu^{-1}(c)$ for some $c\in Z(\g^*)$.

\end{proposition}

\noindent \pf Let $v,w\in \g$ and $X_v, Y_w$ be the induced vector
fields on the manifold $\tsn$, i.e. the flow of $X_v$ is $m\in
\tsn\to e^{vt}.\ m$ and similarly, the flow of $Y_w$ is $m\in
\tsn\to e^{wt}. \ m$. Then:
\begin{equation}\label{iso}
\omega_{St}(X_v,Y_w)=(X_v\lrcorner\omega_{St})(Y_w)=d\mu_v(Y_w)={\mathcal
L}_{Y_w}(\mu_v)
\end{equation}
 \noindent Since ${\mathcal O}$ is a $G$-orbit, the
vector fields $X_v$ for any $v\in \g$ generate all the tangent
space of the orbit. So, ${\mathcal O}$ is isotropic if and only if
$\omega(X_v,Y_w)=0$ for any $v,w \in \g$. Equation $\ref{iso}$
implies that ${\mathcal L}_{Y_w} \mu=0$ on ${\mathcal O}, \forall
w\in \g$ which means that the moment map $\mu$ is constant on
${\mathcal O}$, i.e. ${\mathcal O}\subseteq \mu^{-1}(c)$, for some
$c\in \g^*$. Since ${\mathcal O}$ is $G$-invariant, $c\in
Z(\g^*)$.   $\Box$
\smallskip

 \noindent{\it Remark:} The result ensures that all the
isotropic $G$-orbits, in particular also the Lagrangian orbits are
contained in the level sets of the moment map. Since special
Lagrangian submanifolds are in particular Lagrangian, we have the
following corollary:

\begin{cor}
\label{levset} If $L$ is a connected special Lagrangian
submanifold in $\tsn$ with symmetry group $G\subseteq SO(n+1)$,
then $L\subseteq\mu^{-1}(c)$ for some $c\in Z(\g^*)$, where
$\mu:\tsn\to \g^*$ is the moment map of the action of $G$.

\end{cor}

\section{Homogeneous Special Lagrangian Submanifolds of
$\ts3$}\label{homog}

In this section we will show that the only homogeneous special
Lagrangian $3$-fold in $\ts3$ is the zero section $S^3=\{(x,0)\in
\R^4\times \R^4| |x|=1\}$.

We start by looking at the subgroups of $SO(4)$ that act on the
complex quadric $\ts3$ with generic orbits of dimension 3. Let $U$
be the group of unit quaternions and $\Phi$ the well-known $2:1$
homomorphism $\Phi\colon U\times U\to
 SO(4)$ given by $\Phi(u_1,u_2)(x)=u_1x \bar
 u_2$. It is easy to see that ${\mathfrak so}(4) \cong
{\mathfrak so}(3)_1 \oplus {\mathfrak so}(3)_2$, where $
{\mathfrak so}(3)_1 $ and ${\mathfrak so}(3)_2 $ are two different
copies of ${\mathfrak so}(3)$ whose intersection is the zero
vector. By looking at the subgroups of $SO(4)$ of dimension $\geq
3$ (see also \cite{io}), one can see that the only connected
subgroups of $SO(4)$ that act on $Q$ with generic orbits of
dimension 3 are:

\ni 1. The full group SO(4) with Lie algebra ${\mathfrak so}(4)$
and whose infinitesimal generators are given by:
\begin{equation}\label{genso4}
\l\{\l(\begin{smallmatrix}0&-1&0&0\\1&0&0&0\\0&0&0&0\\0&0&0&0\end{smallmatrix}\r),
\l(\begin{smallmatrix}0&0&-1&0\\0&0&0&0\\1&0&0&0\\0&0&0&0\end{smallmatrix}\r),
\l(\begin{smallmatrix}0&0&0&-1\\0&0&0&0\\0&0&0&0\\1&0&0&0\end{smallmatrix}\r),
\l(\begin{smallmatrix}0&0&0&0\\0&0&1&0\\0&-1&0&0\\0&0&0&0\end{smallmatrix}\r),
\l(\begin{smallmatrix}0&0&0&0\\0&0&0&1\\0&0&0&0\\0&-1&0&0\end{smallmatrix}\r),
\l(\begin{smallmatrix}0&0&0&0\\0&0&0&0\\0&0&0&1\\0&0&-1&0\end{smallmatrix}\r)\r\}
\end{equation}

\ni The generic orbit of the action on $Q$ is an $S^3$.

\ni 2. The subgroup $\widetilde {SO(3)}$, with Lie algebra
${\mathfrak so}(3)_1$ (or ${\mathfrak so}(3)_2$), whose
infinitesimal generators are given by:
\begin{equation}\label{genso3tilde}
\l\{\l(\begin{smallmatrix}0&-1&0&0\\1&0&0&0\\0&0&0&-1\\0&0&1&0\end{smallmatrix}\r),
\l(\begin{smallmatrix}0&0&-1&0\\0&0&0&1\\1&0&0&0\\0&-1&0&0\end{smallmatrix}\r),
\l(\begin{smallmatrix}0&0&0&-1\\0&0&-1&0\\0&1&0&0\\1&0&0&0\end{smallmatrix}\r)\r\}
\end{equation}

\ni 3. The subgroup $S^1\times SO(3)$, with Lie algebra
${\mathfrak so}(2)_1\oplus {\mathfrak so}(3)_2$ (or ${\mathfrak
so}(3)_1\oplus{\mathfrak so}(2)_1$), whose infinitesimal
generators are given by:
\begin{equation}\label{gens1so3}
\l\{\l(\begin{smallmatrix}0&-1&0&0\\1&0&0&0\\0&0&0&0\\0&0&0&0\end{smallmatrix}\r),
  \l(\begin{smallmatrix}0&-1&0&0\\1&0&0&0\\0&0&0&-1\\0&0&1&0\end{smallmatrix}\r),
\l(\begin{smallmatrix}0&0&-1&0\\0&0&0&1\\1&0&0&0\\0&-1&0&0\end{smallmatrix}\r),
\l(\begin{smallmatrix}0&0&0&-1\\0&0&-1&0\\0&1&0&0\\1&0&0&0\end{smallmatrix}\r)\r\}
\end{equation}

\ni We prove the following result:

\begin{proposition}
Every homogeneous special Lagrangian $3$-fold in $T^*S^3$ is
conjugate under the action of $SO(4)$ to the zero section
$S^3\subset T^*S^3$.
\end{proposition}
\pf Let $L$ be a homogeneous special Lagrangian submanifold and
$G\subset SO(4)$ its symmetry group. Then $G$ is one of the three
subgroups of $SO(4)$ described above. We consider each of the
cases.
\smallskip

\ni 1. $G=SO(4)$. From equation (\ref{action}), the moment map of
the $SO(4)$-action is given by $\mu:Q\to {{\mathfrak so}(4)}^*$
with:
$$\mu(z_0,z_1,z_2,z_3)=u'(|z|^2)(\Im(z_0\bar{z_1}),\Im(z_1\bar{z_2}),\Im(z_2\bar{z_3}),\Im(z_3\bar{z_0}),
\Im(z_1\bar{z_3}),\Im(z_2\bar{z_0}))$$

 \ni Since $Z({\mathfrak
so}(4)^*)=\{0\}$, it follows from Corollary \ref{levset} that any
$SO(4)$-invariant special Lagrangian $3$-fold in $Q^3$ lies in the
level set $\mu^{-1}(0)$. By applying an appropriate rotation with
an element of $SO(4)\,$ one can assume that
$x=(x_0,x_1,x_2,x_3)=(\cot t,\sin t,0,0),\ t\in[0,\pi)$. Now,
since the special Lagrangian has to be in the zero level set of
the moment map above, we have
$$\Im(z_0\bar{z_1})=\Im(z_1\bar{z_2})=\Im(z_2\bar{z_3})=\Im(z_3\bar{z_0})=\Im(z_1\bar{z_3})=\Im(z_2\bar{z_0})$$
Using the diffeomorphism (\ref{h}) between the complex quadric
$Q^3$ and $\ts3$, we see that $\xi$ has to be of the form
$\xi=\ro(-\sin t ,\cos t ,0,0)$, where $\ro=|\xi|$ and it has to
also satisfy $\Im(z_0\bar{z_1})=0$, i.e. $\rho=0$. Therefore, $L$
is the zero section of the cotangent bundle.
\smallskip

\ni 2. $G=\widetilde{SO(3)}$. Note that the infinitesimal
generator of $G$ are some linear combinations of the infinitesimal
generator of $SO(4)$. From equation (\ref{action}), the moment map
of the $\widetilde{SO(3)}$-action is given by $\mu:Q\to {\mathfrak
so}(3)^*$ with:
$$\mu(z_0,z_1,z_2,z_3)=u'(|z|^2)(\Im(z_0\bar{z_1}+z_2\bar{z_3}),\Im(z_0\bar{z_2}+z_3\bar{z_1}),
\Im(z_0\bar{z_3}+z_1\bar{z_2}))$$ Since $Z({\mathfrak
so}(3)^*)=\{0\}$, it follows from Corollary \ref{levset} that any
$\widetilde {SO(3)}$-invariant special Lagrangian $3$-fold in
$Q^3$ lies in the level set $\mu^{-1}(0)$. Since the action of
$\widetilde{SO(3)}$ which is given by left multiplication by unit
quaternions is transitive on the unit sphere $S^3$, one can assume
that $x=(x_0,x_1,x_2,x_3)=(\cot t,\sin t,0,0),\ t\in[0,\pi)$. Now,
since the special Lagrangian has to be in the zero level set of
the moment map above,we have
$$\Im(z_0\bar{z_1}+z_2\bar{z_3})=\Im(z_0\bar{z_2}+z_3\bar{z_1})=
\Im(z_0\bar{z_3}+z_1\bar{z_2})$$
 Using this and the diffeomorphism (\ref{h}) yields:
\begin{align}\label{xi}
x_0\xi_1&=x_1\xi_0\notag\\
x_0\xi_2&=x_1\xi_3\\
x_0\xi_3&=-x_1\xi_2\notag
\end{align}

\ni The last two equations yield: $\xi_2=\xi_3=0$. Therefore,
$\xi=\ro(-\sin t ,\cos t ,0,0)$, where $\ro=|\xi|$. First equation
in (\ref{xi}) now implies: $\ro\cos^2 t=-\ro\sin^2 t$, i.e.
$\ro=0$.

\smallskip

\ni 3. $G=S^1\times {SO(3)}$. An argument similar to the previous
cases will also work and this is left to the reader to check.
$\Box$

\bigskip

\section{Cohomogeneity one Special Lagrangian Submanifolds  in
$\ts3$} \ni The next most symmetric case is when the symmetry
group of the special Lagrangian submanifold acts with
cohomogeneity one. In this case, the differential equation of a
special Lagrangian simplifies and we can hope to find examples by
solving an O.D.E. The idea is to find subgroups $G$ of $SO(4)$
whose generic orbits in $\ts3$ are of dimension $2$. In order for
the generic orbit to be dimension $2$, one must have that dim
$G\geq 2$. Then, the strategy is to find an extra direction in
which the submanifold will become special Lagrangian. We will do
this in this section.
\smallskip

\ni The subgroups of $SO(4)$ that act on $Q$ with orbits of
dimension $2$ are the subgroup $SO(3)$ which leaves a direction
invariant and the maximal torus $T^2$. We start out with the
maximal torus case, which is the most interesting case and our
main result.

\subsection{The $T^2$-case}

Let $G$ be the maximal torus $T^2$ of $SO(4)$, described as:
$$\l\{\begin{pmatrix} \cos\theta_1 & -\sin\theta_1 & 0 &0 \\ \sin\theta_1 &\cos\theta_1 & 0 & 0\\ 0& 0 &
\cos\theta_2 & -\sin\theta_2 \\ 0& 0& \sin\theta_2 & \cos\theta_2
\\\end{pmatrix}, \ \theta_1,\theta_2 \in [0,2\pi)\r\}$$ The
following result gives all the special Lagrangian $3$-folds of
$T^*S^3$ invariant under the action of $T^2$.
\begin{theorem}
The special Lagrangian submanifolds in $\ts3 = Q = \{ z \in \C^4
\big| \sum_{i=0}^{3}z_i^2 = 1 \} $ with the Calabi-Yau metric,
which are invariant under the action of the maximal torus $T^2$ of
$SO(4)$ are given by the equations:

\begin{align}
& u'(|z|^2)\,\Im(z_0\bar{z_1}) = c_1 \notag\\
& u'(|z|^2)\,\Im(z_2\bar{z_3}) = c_2 \label{eqt2}\\
& \Im (z_0^2 + z_1^2)= c_3\notag
\end{align}

\ni where $u$ is given by (\ref{diff3}) and $c_1,c_2$ and $c_3$
are any real constants. \label{t}
\end{theorem}

\medskip
\noindent

\ni \pf

\ni Using the two infinitesimal generators of the $T^2$-action:

$$B_1=\l(\begin{smallmatrix}0&-1&0&0 \\
1& 0&0&0\\0&0&0&0\\0&0&0&0\end{smallmatrix}\r), \ B_2=\l(\begin{smallmatrix}0&0&0&0 \\
0& 0&0&0\\0&0&0&-1\\0&0&1&0\end{smallmatrix}\r)$$

\ni and the expression $(\ref{action})$, the moment map for the
$T^2$-action on $Q^3$ is given by:

$$\mu:Q^3\to ({\mathfrak t}^2)^*, \ \mu(z_0,z_1,z_2,z_3)=
u'(|z|^2)(\Im(z_0\bar{z_1}),\Im(z_2\bar{z_3}))$$ Since
$({\mathfrak t}^2)^*=\R^2$, it follows from Corollary
$\ref{levset}$ that any $T^2$-invariant special Lagrangian
$3$-fold $L$ in $Q^3$ lies in a level set $\mu^{-1}(c)$, where
$c=(c_1,c_2)\in \R^2$. The first two equations enforce this
condition and ensures that the submanifold is Lagrangian. One can
also check directly that $\omega_{St}\mid_L=0$ using expression
$(\ref{St2})$. In order to impose the special Lagrangian condition
at a given point $z$, we compute $\Omega_{St}$ on the three
tangent vectors $Y_1=B_1\,z, Y_2=B_2\,z$ and $Y_3= \dot{z}$:

\begin{align*}
 \Omega_{St}(Y_1,Y_2,Y_3)&=(dz_0\wedge dz_1\wedge dz_2\wedge d z_3)(z,Y_1,Y_2,Y_3)= \l| \begin{matrix}
z_0 & - z_1  &  0    & \dot{z_0} \\
z_1 &   z_0  &  0    & \dot{z_1} \\
z_2 &     0  &  -z_3 & \dot{z_2} \\
z_3 &     0  &   z_2 & \dot{z_3}
\end{matrix} \r| \\
&= (z_0^2+z_1^2)(z_2\dot{z_2} + z_3\dot{z_3})-
(z_2^2+z_3^2)(z_0\dot{z_0} + z_1\dot{z_1})
\end{align*}

\smallskip
\ni Now using $z_0^2 + z_1^2 + z_2^2 + z_3^2 = 1, \,\, \Im(z_0^2
+z_1^2) = - \Im(z_2^2+z_3^2)$ and $\Im(z_2\dot{z_2} +
z_3\dot{z_3})= -\Im(z_0\dot{z_0} + z_1\dot{z_1})$ we finally
obtain $\Im\Omega_{St}(Y_1,Y_2,Y_3)= \Im(z_0\dot{z_0} +
z_1\dot{z_1})$, from which the last equation follows. $\Box$

\bigskip

\ni {\it Remark:} Equations (\ref{eqt2}) are obviously
$T^2$-invariant and linearly independent. Therefore, the above
family of $T^2$-invariant special Lagrangian $3$-folds foliate the
cotangent bundle of the sphere, including the zero section. The
generic orbit is $T^2 \times \R$ where $T^2$ is an orbit of the
maximal torus in the isometry group.

\bigskip

One can also view the special Lagrangian $3$-folds above as being
obtained by rotating a curve in $T^*S^3=\{(x,\xi)\in\R^{4}\times
\R^{4}, \
 |x|=1, \ x\cdot \xi=0\}$ by the torus action. To
see this, let $\gamma(t)=(x(t),\xi(t))\in T^*S^3$ be a curve in
the complex quadric. By applying an appropriate rotation with an
element of $T^2$, we can assume that $x(t)=\l(\begin{smallmatrix}\cos t\\
0\\\sin t\\ 0\end{smallmatrix}\r), t\in [0,\pi)$ (since $|x|=1$).
Denote the length of the vector $\xi=\l(\begin{smallmatrix}\xi_0\\
\xi_1\\ \xi_2\\ \xi_3
\end{smallmatrix}\r)$ by $\ro=|\xi|\geq 0$. Let
$\ro_0=\xi_0^2+\xi_2^2$ and $\ro_1=\xi_1^2+\xi_3^2$. Since
$\ro_0^2+\ro_1^2=\ro^2$, we let:
\begin{align*}
\ro_0&=\ro\cos \varphi\\
 \ro_1&=\ro\sin \varphi
\end{align*}
Since $x\cdot \xi=0$, we can parameterize the vector as
$\xi(t)=\l(\begin{smallmatrix} -\ro_0\sin t\\  \rho_1\cos \psi\\
\rho_0\cos t\\  \rho_1 \sin \psi\end{smallmatrix}\r)=\ro\l(\begin{smallmatrix} -\cos \varphi \sin t\\  \sin \varphi \cos \psi\\
\cos \varphi \cos t\\  \sin \varphi \sin
\psi\end{smallmatrix}\r)$. Using the diffeomorphism $h$ given by
relation (\ref{h}), one gets that any point on the quadric is
conjugate under the $T^2$-action to a point of the form:
$$z=\l(\begin{matrix} \cos t\cosh \ro-i\sinh\ro\cos \varphi\sin t\\  i\sinh \ro\sin \varphi \cos \psi\\
 \sin t\cosh \ro+i\sinh\ro\cos \varphi\cos t\\  i\sinh \ro\sin \varphi \sin
\psi\end{matrix}\r)\in Q$$ Note that $|z|^2=\cosh(2\ro)$.

\medskip
\ni In fact the whole quadric $Q^6$ can be parametrized as:

$$ \l(\begin{matrix} \cos\theta_1\cos t\cosh \ro - i\,(\cos\theta_1\sinh\ro\cos \varphi\sin t + \sin \theta_1\sinh \ro\sin \varphi \cos \psi) \\ \sin\theta_1\cos t\cosh \ro + i\,(\cos\theta_1 \sinh \ro\sin \varphi \cos \psi - \sin\theta_1\sinh \ro\cos\varphi \sin t) \\
 \cos\theta_2 \sin t\cosh \ro + i\,(\cos\theta_2\sinh\ro\cos \varphi\cos t - \sin\theta_2\sinh \ro\sin \varphi \sin\psi) \\  \sin\theta_2 \sin t\cosh \ro + i\,(\sin\theta_2\sinh\ro\cos \varphi\cos t + \cos\theta_2\sinh \ro\sin \varphi \sin\psi) \\
 \end{matrix}\r)$$

\smallskip
\ni where $t,\theta_1,\theta_2,\varphi,\psi \in S^1$ and $\ro \geq
0$.
 \bigskip
Equations (\ref{eqt2}) become:
\begin{align}
&u'(\cosh(2\ro))\sinh(2\ro)\cos t\sin \varphi \cos \psi=c_1\notag\\
&u'(\cosh(2\ro))\sinh(2\ro)\sin t\sin \varphi \sin \psi=c_2\label{t2param}\\
&\sinh(2\ro)\sin(2t)\cos \varphi=c_3\notag
\end{align}
These equations are independent of the torus parameters $\theta_1,
\theta_2$ and describe a curve in the parameter space
$(\ro,t,\varphi,\psi)$ which under the $T^2$-action on $Q$ gives a
family of special Lagrangian $3$-folds $L_c$ in the complex
quadric, where $c=(c_1,c_2,c_3)$. These special Lagrangians have
the topology of $T^2\times \R$ and foliate the ambient space $Q$.

\bigskip

{\it Remark:} Another method to find the $T^2$-invariant special
Lagrangians is to start with a curve in the
$(\ro,t,\varphi,\psi)$-space which lies in the level set of the
moment map and write an O.D.E. to impose the special Lagrangian
condition. This O.D.E. is in general difficult to solve explicitly
and the implicit method of proving theorem \ref{t} is much
simpler.

\bigskip
{\bf Asymptotic Behaviour}

\smallskip

\ni We now study the asymptotic behaviour of this family of
special Lagrangian, i.e. the limiting behavior of the family as
$\ro=|\xi|\to\infty$. As we have seen previously, the cotangent
bundle $T^*S^3=\{(x,\xi)\in\R^{4}\times \R^{4}, \ |x|=1, \ x\cdot
\xi=0\}$ approaches the conifold $Q_0$ asymptotically as $|\xi|\to
\infty$.

\ni Notice that as $\ro\to\infty$, $u'(\cosh(2\ro))\to\infty$ from
relation (\ref{diffu}), so in the limit, equations (\ref{t2param})
become one of the following cases:
\begin{align}
&a) \sin \varphi=0, \ \sin t=0\notag\\
&b)  \sin \varphi=0, \ \cos t=0\notag\\
&c)  \cos t=0,\ \sin \psi=0\label{cases}\\
&d) \cos \psi=0,\  \sin t=0\notag
\end{align}
\smallskip

\ni We will study each case separately.\\
a) $ \sin \varphi=0, \ \sin t=0\Rightarrow z=\l(\begin{smallmatrix}\cosh \ro\\
0\\ i\sinh\ro\\ 0 \end{smallmatrix}\r)$. The unit vector
$\frac{z}{|z|}$ is:
$$\frac{z}{|z|}= \frac{1}{\sqrt{\cosh(2\ro)}} \l(\begin{matrix}\cosh \ro\\
0\\ i\,\sinh\ro\\ 0
\end{matrix}\r)\to \frac{1}{\sqrt{2}}\l(\begin{matrix} 1\\
0\\ i\\ 0
\end{matrix}\r)\in Q_0, \ \mbox{ as } \ro\to \infty $$
Applying the $T^2$-action in the limit, one gets a surface
$\Sigma$ diffeomorphic to $T^2$ and $\Sigma$ is a submanifold of
the conifold $Q_0$:
$$\Sigma=\l\{\frac{1}{\sqrt{2}}\l(\begin{matrix}\cos \theta_1\\
\sin \theta_1\\  i \, \cos \theta_2\\
i\,\sin \theta_2\end{matrix}\r), \theta_1,\theta_2\in[0,2\pi)
\r\}\subset Q_0$$

\ni We will show that the cone on $\Sigma, \ C(\Sigma)=\l\{sz|z\in
\Sigma, \ s\in\R\r\}$, is special Lagrangian in the conifold
$Q_0$, endowed with the Ricci-flat metric found by Candelas and de
la Ossa\cite{co}.

\ni We first show that the cone $C(\Sigma)$ is Lagrangian, i.e.
$\omega_{cone}\mid_{C(\Sigma)}=0$. The moment map of the
$T^2$-action on the cone is: $$\mu_0:Q_0\to \R^2, \
\mu_0(z)=u_{cone}'(|z|^2)(\Im(z_0\bar{z_1}),\Im(z_2\bar{z_3}))$$
where $u_{cone}(r^2)$ is the potential function for the conifold
given in section \ref{sect}. Since the cone on $\Sigma$ is seen to
lie in $\mu_0^{-1}(0,0)$, it follows from Prop. \ref{lagr} that it
is Lagrangian.

Next we show that the cone is special Lagrangian, i.e.
$\Im\Omega_{cone}\mid_{C(\Sigma)}\equiv0$. For this we compute
$\Omega_{cone}$ on three tangent vectors $Y_1,Y_2,Y_3$ to the cone
$C(\Sigma)$. One of them is the position vector and the other two
vectors are the derivatives with respect to the parameters
$\theta_1$ and $\theta_2$. The unit vector normal to the cone is
given by
$$w=\bar{z}=\frac{1}{\sqrt{2}}\l(\begin{matrix} \cos \theta_1\\
 \sin \theta_1\\ -i \, \cos \theta_2 \\
-i \, \sin \theta_2  \end{matrix}\r)$$ and this is the vector we
will use to compute $\Omega_{cone}$ as follows:

\begin{align*}
 \Im\Omega_{cone}(Y_1,Y_2,Y_3)&=\Im(dz_0\wedge dz_1\wedge dz_2\wedge d z_3)(\bar{z},Y_1,Y_2,Y_3) \\
&= \Im\frac{1}{2s^2}\l| \begin{matrix}
\frac{\cos \theta_1}{\sqrt{2}}& \frac{\cos \theta_1}{\sqrt{2}} &  -s\frac{\sin \theta_1}{\sqrt{2}}    & 0\\
 \frac{\sin \theta_1}{\sqrt{2}} &  \frac{\sin \theta_1}{\sqrt{2}} &  s\frac{\cos\theta_1}{\sqrt{2}}   & 0 \\
 -\frac{i\cos \theta_2}{\sqrt{2}} &   \frac{i\cos \theta_2}{\sqrt{2}}   &  0&
 -s \frac{i\sin \theta_2}{\sqrt{2}}  \\
  -\frac{i\sin \theta_2}{\sqrt{2}}&     \frac{i\sin \theta_2}{\sqrt{2}}  &
   0 &  s\frac{i\cos \theta_2}{\sqrt{2}}\end{matrix} \r| =0
\end{align*}
Since $\Im\Omega_{cone}\mid_{C(\Sigma)}\equiv0$, the cone on
$\Sigma$ is special Lagrangian.

{\it Remark:} The above analysis shows that in this case, the
$T^2$-invariant special Lagrangian in the quadric goes
asymptotically to a special Lagrangian cone in the conifold.
\smallskip

\ni b) $ \sin \varphi=0, \ \cos t=0$. In this case, in the limit
as $\ro\to \infty$, the unit vector $\frac{z}{|z|}$ goes to:
$$ \frac{1}{\sqrt{2}}\l(\begin{matrix} -i\\
0\\ 1 \\ 0
\end{matrix}\r)\in Q_0 $$
and the conclusion is the same as in case a). The unit vector
$\frac{z}{|z|}$ in this case is the unit vector from case a)
rotated by $J$ and so the the limit is the same special Lagrangian
cone.

\smallskip
\ni c) $ \cos t=0, \ \sin \psi=0\Rightarrow z=\l(\begin{smallmatrix}-i\sinh\ro\cos \varphi\\
 i\sin \varphi\sinh\ro\\ \cosh \ro\\0 \end{smallmatrix}\r)$. The unit
vector $\frac{z}{|z|}$ goes in the limit to:
$$\frac{1}{\sqrt{2}}\l(\begin{matrix}-i\cos \varphi\\
i\sin \varphi\\ 1\\ 0
\end{matrix}\r)\in Q_0$$
Applying the $T^2$-action in the limit, one gets a surface
$\Sigma$ diffeomorphic to $T^2$ and $\Sigma$ is a submanifold of
the conifold $Q_0$:
$$\Sigma=\l\{ \frac{1}{\sqrt{2}}\l( \begin{matrix} i \, \cos \theta_1 \\
i \, \sin\theta_1 \\ \cos \theta_2 \\
\sin \theta_2 \end{matrix}\r), \theta_1,\theta_2 \in[0,2\pi)
\r\}\subset Q_0$$ Same argument as in case a) shows that the cone
on $\Sigma$ is special Lagrangian in the conifold $Q_0$.

\smallskip

 \ni d) $ \cos \psi=0, \ \sin t=0$. This case yields the same special Lagrangian cone as in case c) and it is
 left to the reader.

\smallskip

\ni {\it Remark 1:} The special Lagrangian family in the complex
quadric is asymptotic to a special Lagrangian cone on flat tori in
the conifold. For a fixed $c=(c_1,c_2,c_3)$, $L_c$ has two
components, each of them asymptotic to the special Lagrangian cone
in the conifold.

\smallskip

\ni {\it Remark 2:} When $\rho=0$, equations $(\ref{t2param})$ are
identically satisfied and the solution is the zero section $S^3$
of the cotangent bundle, which is well-known to be special
Lagrangian.

\smallskip

\ni {\it Remark 3:}  If we set $c_1=c_2=0$ in equation
(\ref{t2param}), we obtain the equation $$\sin (2t)\sinh(2\ro)=c$$
in the $(t,\ro)$-plane, $t\in [0,\pi)$. The special Lagrangian in
this case lies in the zero-level set of the moment map. The phase
portrait of the equation in this special case is shown in Figure
1. We can see that as $t\to\frac{\pi}{2}$, $|\rho|\to \infty$. The
special Lagrangians are obtained by rotating these curves by a
$T^2$-action.
\smallskip
\begin{figure}
\begin{center}
\includegraphics{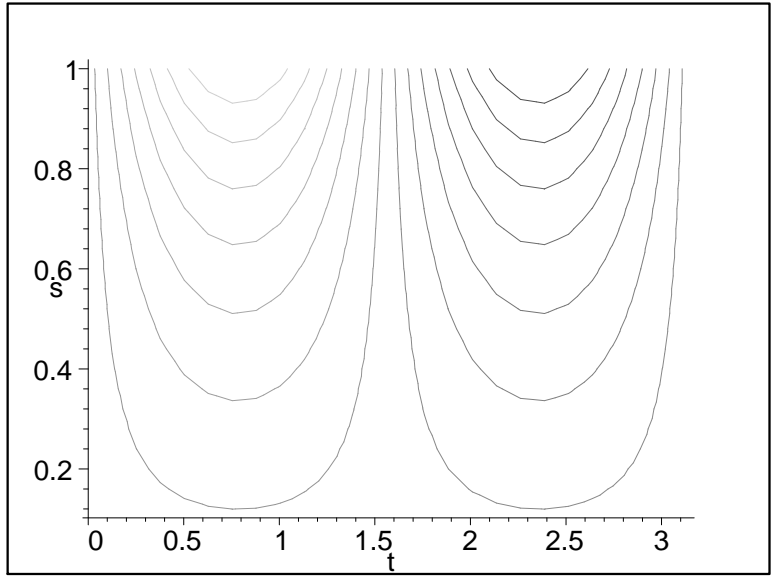}
\end{center}
\caption{\label{Figure 1} Phase Portrait for the equation $\sin
(2t)\sinh(2s)=c$ }
\end{figure}

\subsection{The $SO(3)$-case}

Let $G$ be the subgroup $SO(3,\R)$ of $SO(4,\C)$, where $SO(3,\R)$
sits in $SO(4,\R)$ as matrices of the form
$$\l\{\begin{pmatrix}1&0 \\ 0& A\end{pmatrix},\ A\in SO(3)\r\}$$
and $SO(4,\R)$ is embedded in $SO(4,\C)$ as $4 \times 4$ real
matrices. Note that this $SO(3)$ is embedded differently in
$SO(4)$ than the subgroup $\widetilde{SO(3)}$ from section
\ref{homog}.

\begin{theorem}\label{thm1}
Let $\ts3=\{(x,\xi)\in\R^{4}\times \R^{4}, \ |x|=1, \ x\cdot
\xi=0\}$. Then, $$L_c=\{(gx,g\xi)\mid (x,\xi)\in\ts3,\ g\in
SO(3),\ x=(\cos t,\sin t,0,0), \
2|\xi|-\cos(2t)\sinh(2|\xi|)=c\}$$ is the only $SO(3)$-invariant
special Lagrangian $3$-fold of $\ts3$.
\end{theorem}

\ni \pf Let $z=(z_0,z_1,z_2,z_3)\in Q^3$, i.e.
$z_0^2+z_1^2+z_2^2+z_3^2=1$. Using expression (\ref{action}), the
moment map of the $SO(3)$-action on $Q^3$ is given by:
$$\mu:Q^3\to {\mathfrak so}(3)^*, \ \mu(z_0,z_1,z_2,z_3)=
 u'(|z|^2)(\Im(z_1\bar{z_2}),\Im(z_2\bar{z_3}),\Im(z_3\bar{z_1}))$$
To see how we obtain this, let
$$ A_1=\l(\begin{smallmatrix}0&0&0&0 \\
0& 0&0&0\\0&0&0&-1\\0&0&1&0\end{smallmatrix}\r), \,
A_2=\l(\begin{smallmatrix}0&0&0&0 \\
0& 0&0&1\\0&0&0&0\\0&-1&0&0\end{smallmatrix}\r), \,
A_3=\l(\begin{smallmatrix}0&0&0&0 \\
0& 0&-1&0\\0&1&0&0\\0&0&0&0\end{smallmatrix}\r)$$ be the
infinitesimal generators of the $SO(3)$-action. Equation
$(\ref{action})$ implies that:
\begin{align*}
\mu_{A_3}(z)&=\frac{1}{2}u'(|z|^2)<A_3z,iz>=\frac{1}{2}u'(|z|^2)\l<\begin{pmatrix}
0\\-z_2\\z_1\\0
\end{pmatrix},\begin{pmatrix}
iz_0\\iz_1\\iz_2\\iz_3\end{pmatrix}\r>=\\
&=\frac{1}{2}u'(|z|^2)(i\bar{z_1}z_2- i z_1\bar{z_2})=
u'(|z|^2)\Im(z_1\bar{z_2})
\end{align*}
Similarly, $\mu_{A_2}(z)= u'(|z|^2)\Im(z_3\bar{z_1})$,
$\mu_{A_1}(z)= u'(|z|^2)\Im(z_2\bar{z_3})$ and we get the moment
map of the $SO(3)$-action on $Q^3$.

\noindent Since $Z({\mathfrak so}(3)^*)=\{0\}$, it follows from
Corollary \ref{levset} that any $SO(3)$-invariant special
Lagrangian $3$-fold in $Q^3$ lies in the level set $\mu^{-1}(0)$.
By applying an appropriate rotation with an element of $SO(3)$, we
can assume that $x=(x_0,x_1,x_2,x_3)=(\cot t,\sin t,0,0),\
t\in[0,\pi)$. Now, since the special Lagrangian has to lie in the
zero level set of the moment map,we have
$\Im(z_1\bar{z_2})=\Im(z_2\bar{z_3})=\Im(z_3\bar{z_1})=0$ and
hence the level set $\mu^{-1}(0)$ is given by:
 $$\mu^{-1}(0)= \{ \,(gx,g\xi) | g \in SO(3), \, x(t)=(\cos
t,\sin t ,0,0), \, \xi(t)=\ro(-\sin t ,\cos t ,0,0), \, t \in [0,
\pi), \, \ro \geq 0 \, \}$$

\ni or equivalently, using the identification of $\ts3$ with $Q$,
$$ \l\{g.z(t,\rho) \,: \,\, g\in SO(3), \ z(t,\ro)=
\begin{pmatrix} \cos t \cosh\ro-i \sin t\sinh \ro \\
\sin t\cosh \ro+i \cos t \sinh \ro \\ 0
\\0 \end{pmatrix}
= \begin{pmatrix}\cos (\tau)\\
\sin (\tau) \\ 0 \\0 \end{pmatrix} ,\, \ro \geq 0, \, t\in[0,\pi)
\r\}$$ \ni where $\tau=t+i\ro$ lies in the $[0,\pi)$ vertical
strip of the complex plane.

\medskip
\noindent We now look for curves $\gamma(s)$ in the $\tau =
(t,\ro)$ plane which after applying the $SO(3)$-action give rise
to special Lagrangians of the form $L = g.\gamma(s)$.
If $p = \gamma(s) = \begin{pmatrix}\cos (\tau(s))\\
\sin (\tau(s)) \\ 0
\\0\end{pmatrix}$ is a point on the level set $\mu^{-1}(0)$,
then  $A_3\,p=\begin{pmatrix}0\\ 0\\ \sin(\tau)
\\0\end{pmatrix}$,
$A_2\,p=\begin{pmatrix} 0 \\
 \\ 0 \\ -\sin(\tau) \end{pmatrix}$ and  $A_1\, p= 0 \,$,
where $A_1,A_2,A_3$ are the infinitesimal generators of
${\mathfrak so}(3)$ as defined above.

\medskip
\ni The tangent plane at $p$ to $L$ is spanned by the vectors
$<X_1=A_1\,p, X_2=A_2\,p, X = \dot{\gamma}(s) >$ at $p$, where $
 \dot{\gamma}(s)= \begin{pmatrix} -\sin (\tau)\, \dot{\tau}\\
\cos (\tau)\, \dot{\tau} \\ 0
\\0\end{pmatrix}$.

\smallskip
\noindent $L$ is invariant under the $SO(3)$-flow and
$\omega_{St}|_{L}=0$, since it lies in the zero level set of the
moment map. Therefore $L$ is Lagrangian. One can also verify this
directly with formula $(\ref{St2})$. Now we will impose the
condition that $L$ is special Lagrangian, i.e. $\Im\Omega_{St}=0$
should hold.

\smallskip
\ni Using equation (\ref{holom}), we compute
$\Omega_{St}(X_1,X_2,X_3)  = \big(d z_0 \wedge d z_1
 \wedge d z_2\wedge d z_3\big)\big(\gamma(s), X_1, X_2, \dot{\gamma}(s) \big)$:

$$ \l| \begin{matrix}
\cos(\tau)& 0 & 0 &-\sin(\tau)\,\dot{\tau}\\
\sin(\tau) & 0 & 0 & \cos(\tau)\,\dot{\tau}\\
0 & \sin(\tau) & 0 & 0 \\
0 & 0 & -\sin(\tau) & 0
\end{matrix} \r|
= - \sin^2(\tau)\, \dot{\tau}
$$

\smallskip
\noindent Integrating, the condition $\Im\Omega_{St}=0$ becomes:
$$\Im(2\tau-\sin(2\tau))=c$$ which is equivalent to
\begin{equation}\label{so3}
2\ro-\cos(2t)\sinh(2\ro)=c
\end{equation}where $c$ is any real
constant.  $\Box$

\medskip

\ni {\it Remark 1 :} Notice that $\rho=0$ is a solution to
equation (\ref{so3}) and the special Lagrangian obtained is the
zero section of $T^*S^3$, which is known to be special Lagrangian.

\smallskip
\ni {\it Remark 2:} The result we obtained in the $n=3$ case can
be generalized to study the $SO(n)$-invariant special Lagrangian
$n$-folds in $\tsn$. We will consider this case in the next
section.

\smallskip
\noindent {\it Remark 3:} This family of $SO(n)$-invariant special
Lagrangian submanifolds was also obtained by Anciaux \cite{an}
using different methods.

\bigskip

The $SO(3)$-invariant special Lagrangian $L$ can also be written
intrinsically as follows. Choose $(z_0,z_1,z_2,z_3)$ coordinates
on the complex quadric $Q^3$. Then $L$ is given by the equations:
\begin{align} \label{eqso3}
&u'(|z|^2)\Im(z_1\bar{z_2})=c_1\notag\\
&u'(|z|^2)\Im(z_2\bar{z_3})=c_2\\
&\Im\l(\arccos(z_0)-z_0\sqrt{1-z_0^2}\r)=c\notag
\end{align}
where $c_1,c_2,c$ are constants.

\bigskip

 \noindent

 \ni {\bf Asymptotic Behaviour:}

\medskip

\ni In what follows, we will study the equation
$2\ro-\cos(2t)\sinh(2\ro)=c$ in the $(t,\ro)$-plane and describe
the asymptotic behaviour of the special Lagrangian $3$-folds
obtained in theorem $\ref{thm1}$. \vskip.6cm

\begin{figure}
\begin{center}
\includegraphics{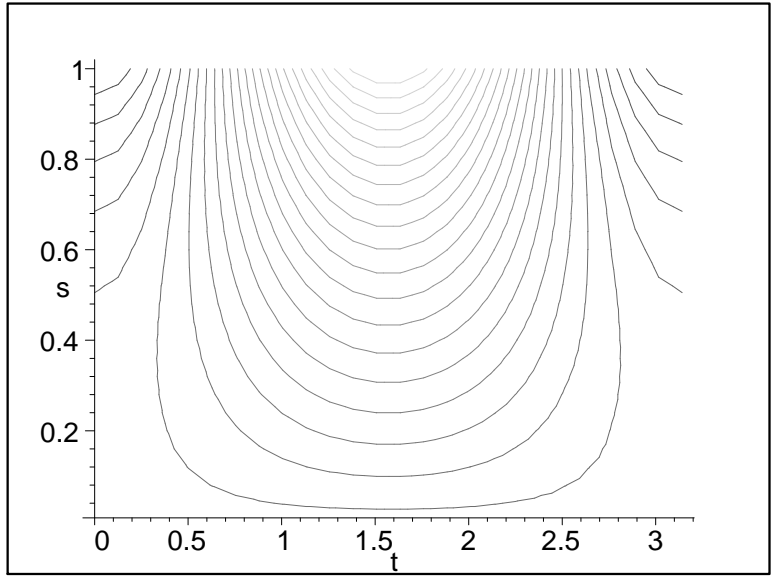}
\end{center}
\caption{\label{Figure 2} Phase Portrait for the equation
$2s-\cos(2t)\sinh(2s)=c$}
\end{figure}


\ni When $\rho\to \infty$, equation (\ref{so3}) becomes $\cos
2t=0$, so $t=\frac{\pi}{4}$. The unit vector $\frac{z}{|z|}$ is:
$$\frac{1}{{\sqrt{ 2 \cosh(2\ro)}}} \l(\begin{matrix} \cosh \ro -i \, \sinh \ro\\
\cosh \ro + i\, \sinh \ro\\
0 \\
0
\end{matrix}\r)
\to \frac{1}{2}\l(\begin{matrix}1-i\\ 1+i\\ 0\\ 0
\end{matrix}\r)\in Q_0, \ \mbox{ as }
\ro \to \infty $$

\ni Applying the $SO(3)$-action in the limit, one gets a surface
$\Sigma$ diffeomorphic to $S^2$ and $\Sigma$ is a submanifold of
the conifold $Q_0$. As before, the cone on $\Sigma$ is special
Lagrangian in the conifold $Q_0$, endowed with the Ricci-flat
metric found by Candelas and de la Ossa\cite{co}. To see this, we
will first show that the cone $C(\Sigma)$ is Lagrangian, i.e.
$\omega_{cone}\mid_{C(\Sigma)}=0$. The moment map of the
$SO(3)$-action on the cone is: $$\mu_0:Q_0\to \R^3, \
\mu_0(z)=u_{cone}'(|z|^2)(\Im(z_1\bar{z_2}),\Im(z_2\bar{z_3}),\Im(z_3\bar{z_1}))$$
where $u_{cone}(r^2)$ is the potential function for the conifold
given in section \ref{sect}. Since the cone on $\Sigma$ is seen to
lie in $\mu_0^{-1}(0,0)$, it follows from Prop. \ref{lagr} that it
is Lagrangian.

Next we show that the cone is special Lagrangian, i.e.
$\Im\Omega_{cone}\mid_{C(\Sigma)}\equiv 0$. For this we compute
$\Omega_{cone}$ on three tangent vectors $Y_1,Y_2,Y_3$ to the cone
$C(\Sigma)$. One of them is the position vector and the other two
vectors are the vectors $A_3z$ and $A_2z$, where $A_2,A_3$ are the
infinitesimal generators as defined above. The unit vector normal
to the cone is given by
$$w=\bar{z}=\frac{1}{2}\l(\begin{matrix}1+i\\
1-i\\ 0\\ 0\end{matrix}\r)$$ and this is the vector we will use to
compute $\Omega_{cone}$ as follows:

\begin{align*}
 \Im\Omega_{cone}(Y_1,Y_2,Y_3)&=\Im(dz_0\wedge dz_1\wedge dz_2\wedge d z_3)(\bar{z},Y_1,Y_2,Y_3) \\
&= \Im\frac{1}{8}\l| \begin{matrix}
1-i &  0   & 0&1+i\\
1+i & 0 & 0  & 1-i \\
0&   1+i   &  0&
 0 \\
0&   0  &
   -1-i&  0\end{matrix}\r| =0
\end{align*}
Since $\Im\Omega_{cone}\mid_{C(\Sigma)}\equiv0$, the cone on
$\Sigma$ is special Lagrangian.

{\it Remark:} The above analysis shows that in this case, the
$SO(3)$-invariant special Lagrangians in the quadric approach
asymptotically a special Lagrangian cone on $S^2$ in the conifold.
This limiting $S^2$ can be explicitly described in coordinates as
follows:

$$z =\frac{1}{2}\l(\begin{matrix}1-i\\
(1+i)\cos \phi \cos \theta\\ (1+i)\cos \phi \sin \theta \\
(1+i)\sin \phi \end{matrix}\r) \qquad \phi \in [-\frac{\pi}{2},
\frac{\pi}{2}], \, \theta \in [0, 2\pi)$$

\bigskip
{\it Remark:} In particular, in the case of Eguchi-Hanson metric
on $Q^2$ ($n=2$), the subgroups $SO(2)$ and $T^1$ of $SO(3)$
coincide and equations (\ref{eqt2}) and (\ref{eqso3}) give the
same special Lagrangian submanifold invariant under $SO(2)$,
embedded in $SO(3)$ as:
$$\l\{\l(\begin{matrix}1&0\\0&A\end{matrix}\r), \ A\in
SO(2)\r\}$$
 In coordinates $z=(z_0,z_1,z_2)$ on $Q$, these special
Lagrangians are given by the equations:
\begin{align*}
&\frac{|z|}{\sqrt{|z|^2+1}}\Im(z_1\bar{z_2})=c_1\\
&\Im (z_0)=c_2
\end{align*}
where $c_1,c_2$ are constants.

\subsection{The general $SO(n)$-case}

One can generalize our method to higher dimensions and recover the
$SO(n)$-invariant family of special Lagrangians that H. Anciaux
obtained in \cite{an}. Computing the imaginary part of the
holomorphic $n$-form $\Omega_{St}$ to check the special Lagrangian
condition yields $$\Im(\sin^{n-1}(\tau) \dot{\tau})=0$$ Now let
$F(\tau)$ be the function $\Im \l( \int_0^\tau sin^{n-1}(\sigma)
d\sigma \r ) $. Combining with the moment map conditions, the
$SO(n)$-invariant special Lagrangians are given by the following
set of equations:

\begin{align} \label{eqson}
&u'(|z|^2)\Im(z_1\bar{z_j})=c_j \qquad  2 \leq  j \leq n \\
&\Im\l(F(\arccos(z_0)) \r)=c \notag
\end{align}
where and $c_j,c$ are constants and the function $u$ satisfies
(\ref{ode}) for the given dimension.

\section{Concluding remarks:}

We conclude the paper with a remark and some comments.
\smallskip

\ni {\it Remark:} There are no special Lagrangian submanifolds in
$\tsn$ which are graphs over the zero section $S^n$.

\ni \pf Let $L=\l\{(x,\xi(x))\in \R^{n+1}\times \R^{n+1}|\ |x|=1,
\ x\cdot \xi=0\r\}$ be a special Lagrangian such that $L$ is a
graph over the zero section $S^3$. Then $L$ is minimal and is
homotopic to the zero section. Since a special Lagrangian is
absolutely area minimizing in its homology class, it follows that
$L$ is the zero section $S^n$.

 The SYZ conjecture \cite{syz} explains mirror symmetry by
analyzing "dual" fibrations of two different Calabi-Yau $3$-folds
by special Lagrangian tori $T^3$. Since in our case we are dealing
with a non-compact Calabi-Yau, we obtain a $T^2\times
\R$-fibration of the deformed conifold instead of a $T^3$. It is
known that the local mirror of $\ts3$ is the ${\mathcal
O}(-1)\oplus {\mathcal O}(-1)$-bundle over $\C P^1$. In our next
paper, we plan to investigate what happens to these special
Lagrangians under the conifold transition to the mirror.


\ni MATHEMATICS DEPARTMENT, MCMASTER UNIVERSITY\\
{\it e-mail address:} ionelm@@math.mcmaster.ca \\

\ni MATHEMATICS DEPARTMENT, MCMASTER UNIVERSITY\\
{\it e-mail address:} minoo@@mcmaster.ca

\end{document}